\documentclass[11pt]{article}
\usepackage{amsmath} 
\usepackage{lmodern}
\usepackage[T1]{fontenc}
\usepackage[english]{babel}
\usepackage[utf8]{inputenc}
\usepackage[a4paper]{geometry}
\usepackage{amssymb}
\usepackage[hypcap]{caption}
\usepackage{amsmath, amsthm, amsfonts, amscd}
\usepackage[nottoc,numbib]{tocbibind}
\usepackage{mathrsfs}
\usepackage{float}
\usepackage{graphicx}
\usepackage[labelformat=empty]{caption}
\usepackage{microtype}
\usepackage{todonotes}
\usepackage{makeidx}
\usepackage{verbatim}
\usepackage{hyperref}
\hypersetup{
colorlinks=false,% hyperlinks will be black
linkbordercolor=red,% hyperlink borders will be red
pdfborderstyle={/S/U/W 1}% border style will be underline of width 1pt
}
\usepackage{mdwlist}
\usepackage{mathabx}
\usepackage{mathtools}
\usepackage{titlesec}

\usepackage[nottoc,numbib]{tocbibind}

\newtheorem{theorem}{\scshape{Theorem}}[section]
\newtheorem{corollary}[theorem]{\scshape{Corollary}}
\newtheorem{proposition}[theorem]{\scshape{Proposition}}
\newtheorem{lemma}[theorem]{\scshape{Lemma}}

\newtheorem*{theorem*}{Theorem}

\newtheorem*{prop*}{Proposition}
\newtheorem*{cor*}{Corollary}

\theoremstyle{definition}
\newtheorem{definition}[theorem]{\scshape{Definition}}

\newtheorem{example}[theorem]{\scshape{Example}}

\newtheorem{remark}[theorem]{\scshape{Remark}}
\newtheorem{problem}[theorem]{\scshape{Problem}}

\def\G{\Gamma}

\def\subset{\subseteq}

\def\={\cong}

\def\mb{\mathbf}
\def\mc{\mathcal}
\def\mbb{\mathbb}

\def\Is{{\it Is}}
\def\ker{{\it ker}}
\def\G_0{G/\Is(\gamma_{3}(G))}

\def\End{\it End}

{}

\title{Diophantine problems in solvable groups}
\author{Albert Garreta\footnote{\noindent University of the Basque Country,  Spain, \emph{garreta.a@gmail.com} (corresponding author).}, Alexei Miasnikov\footnote{\noindent Stevens Institute of Technology, NJ, USA, \emph{amiasnikov@gmail.com}}, and Denis Ovchinnikov\footnote{\noindent Stevens Institute of Technology, NJ, USA, \emph{dovchinn@stevens.edu} \newline
This work was supported by the Mathematical Center in Akademgorodok. 
\newline The first named author was supported by the ERC grant 336983, by the Spanish Government grant MTM2017-86802-P, and by the Basque Government grant IT974-16.\newline
\emph{MSC2010 subject classification}: 20F70, 20F10, 03B25, 03D35, 20F18, 20F16
}}

\date{December, 2019}

\begin{document}

\maketitle

\begin{abstract}
We study the Diophantine problem (decidability of finite systems of equations) in different classes of finitely generated solvable groups (nilpotent, polycyclic, metabelian, free solvable, etc), which satisfy some natural "non-commutativity" conditions. For each group $G$ in one of these classes, we prove that there exists a ring of algebraic integers  $O$ that is interpretable in $G$ by finite systems of equations (e-interpretable), and hence that the Diophantine problem  in $O$ is polynomial time reducible  to the Diophantine problem in $G$.  One of the major open conjectures in number theory states that the Diophantine problem in any such $O$ is undecidable. If true this would imply that the Diophantine problem in any such $G$ is also undecidable.  Furthermore, we show that for many particular groups $G$ as above, the ring $O$ is isomorphic to the ring of integers $\mathbb{Z}$, so the Diophantine problem in $G$ is, indeed, undecidable. This holds, in particular, for free nilpotent or  free solvable  non-abelian groups, as well as  for non-abelian generalized Heisenberg groups and uni-triangular groups $UT(n,\mathbb{Z}), n \geq 3$. Then we apply these results to non-solvable groups that contain non-virtually abelian maximal finitely generated nilpotent subgroups. For instance, we show that the Diophantine problem  is undecidable in the groups $GL(3,\mathbb{Z}), SL(3,\mathbb{Z}), T(3,\mathbb{Z})$.
\end{abstract}

\tableofcontents

\section{Introduction}\label{basic}

We study the Diophantine problem (decidability of finite systems of equations) in different classes of finitely generated solvable groups (nilpotent, polycyclic, metabelian, free solvable, etc), which satisfy some natural "non-commutativity" conditions. For each group $G$ in one of these classes, we prove that there exists a ring of algebraic integers  $O$ that is interpretable in $G$ by finite systems of equations (e-interpretable), and hence that the Diophantine problem  in $O$ is effectively polynomial time reducible to the Diophantine problem in $G$. A famous   conjecture in number theory states that the Diophantine problem in any such $O$ is undecidable, implying by the result above  that the Diophantine problem in any such $G$ would be also undecidable.  In fact, we show that for many particular groups $G$ as above, the ring $O$ is isomorphic to the ring of integers $\mathbb{Z}$, so the Diophantine problem in $G$ is, indeed, undecidable. This holds, in particular, for free nilpotent or  free solvable  non-abelian groups, as well as  for non-abelian generalized Heisenberg groups and uni-triangular groups $UT(n,\mathbb{Z}), n \geq 3$. Then we apply these results to non-solvable groups that contain non-virtually abelian maximal finitely generated nilpotent subgroups. For instance, we show that the Diophantine problem  is undecidable in the groups $GL(3,\mathbb{Z}), SL(3,\mathbb{Z}), T(3,\mathbb{Z})$.

%We apply these results to non-solvable groups that contain non-virtually abelian maximal finitely generated nilpotent subgroups. For instance, we show that Diophantine problem  is undecidable in the groups $GL(3,\mathbb{Z}), SL(3,\mathbb{Z}), T(3,\mathbb{Z})$, and in non-abelian free solvable groups.

The Diophantine problem (also called Hilbert's tenth problem or generalized Hilbert's tenth problem) in a structure $R$, denoted $\mc{D}(R)$, asks whether there exists an algorithm that, given a finite \emph{system} of equations $S$ with coefficients in $R$, determines  if $S$ has a solution in $R$ or not.  The original version of this problem was posed by Hilbert for the ring of integers $\mbb Z$. This was solved in the negative in 1970 by Matiyasevich \cite{mat} building on the work of Davis, Putnam, and Robinson \cite{DPR}. Subsequently the same problem has been studied in a wide variety of rings, most notably in $\mbb Q$ and in rings of algebraic integers $O$ (integral closures of $\mbb Z$ in finite field extensions of $\mbb Q$),
where it remains widely open. A long-standing conjecture  (see, for example, \cite{Denef_conjecture, Phe_Zah}) states that $\mbb Z$ is Diophantine in any such $O$  (and thus  $\mc{D}(O)$ is undecidable). This conjecture  has been verified in some particular cases \cite{Sha_Shla, Shla_book, new}, and it has been shown to be true assuming the Safarevich-Tate conjecture \cite{Mazur2010}. We refer to \cite{Poonen, Phe_Zah, Shla_book} for further information on the Diophantine problem in different rings and fields of number-theoretic flavor. On the other hand, Kharlampovich and the second author showed in \cite{KM_free_algebras}  that the Diophantine problem is undecidable for  in free associative algebras  for any field of coefficients, and in the group algebras (over any field of coefficients) of a wide variety  of torsion-free groups, including
toral relatively hyperbolic groups, right angled Artin groups, commutative
transitive groups, and the fundamental groups of various graph groups. Moreover, they proved in \cite{KMLie} undecidability of the Diophantine problem in free Lie algebras of  rank  at least three with coefficients in  an arbitrary integral domain.   In \cite{GMO_rings} we studied the Diophantine problem in more general rings and algebras (possibly non-associative, non-commutative, and non-unitary), obtaining analogous results to the ones in this paper. Indeed, the present paper may be read as a continuation of \cite{GMO_rings}.

Research on systems of equations and their decidability in groups has a very long history, it goes back to 1912 to the pioneering works of Dehn on the word and conjugacy problems in finitely presented groups.  Within the class of solvable groups, it is known that the word and conjugacy are decidable for many such groups, including finitely generated nilpotent, polycyclic, metabelian, and  free solvable groups. On the other hand, there is also a famous example, due to Kharlampovich,  of a finitely presented solvable group with undecidable word problem \cite{KHAR}. 

The first results on the proper Diophantine problem in groups are due to Romankov. He   showed in \cite{Romankov_free_metabelian, Romankov_eqns_2} that the Diophantine problem is undecidable in any non-abelian free metabelian group and in any non-abelian free nilpotent group of nilpotency class at least $9$. Variations and improvements of these results were obtained subsequently in \cite{Burke, Durnev, Truss}.  Recent work of   Duchin, Liang and Shapiro  \cite{Duchin} shows that $\mc{D}(N)$  is undecidable in any finitely generated  nonabelian free nilpotent group $N$. We refer to a survey \cite{Romankov_survey} for  these and more results on equations in groups. Stepping outside of the realm of systems of equations, Noskov showed in \cite{Noskov}, following the work of Malcev \cite{Malcev_correspondence}, Ershov \cite{Ershov} and Romanovskii \cite{romanovskii},    that the first-order theory of any finitely generated non-virtually abelian solvable group is undecidable. Note,  that  much earlier Ershov \cite{Ershov} proved that any virtually abelian group has decidable elementary theory.  The papers \cite{Rom_eqn} and \cite{Chapuis} contain results of a similar flavour to the ones of this paper:  decidability of the universal theory of a free nilpotent group or a free solvable group of class at least 3 implies decidability of the Diophantine problem in the field of rational numbers $\mbb Q$, a major open problem.

In solvable groups systems of equations are fundamentally different from single equations.  For instance, finite systems of equations are undecidable in the Heisenberg group (i.e.\ the free nilpotent group of nilpotency class $2$ and of rank $2$), while single equations are decidable  \cite{Duchin}. This contrasts with most number theoretic settings, where the two notions are often used interchangeably since  in all integral domains, whose field of fractions are not algebraically closed, every finite system of equations is equivalent to a single equation (see, for example,\cite{Poonen}). Much of the research regarding equations in solvable groups was  focused so far on single equations (see \cite{Romankov_survey}), indeed Romankov's aforementioned results \cite{Romankov_eqns_2, Romankov_free_metabelian}  (and also Truss' \cite{Truss}) actually prove that single equations are undecidable in the corresponding groups $G$. These are stronger results than just undecidability of the Diophantine problem in $G$.   Allowing  arbitrary finite systems of equations is fundamental to our approach. This  makes the whole theory much more robust and brings to the table powerful general methods.  %  

This line of results  changes  drastically outside of the class of solvable groups: the work of Makanin and Razborov \cite{Makanin, Razborov1}  shows that $\mc{D}(F)$ is decidable for any  free group $F$, and it further provides  a  description of the solution sets to  arbitrary systems of equations in $F$ (systems of equations are equivalent to  single equations in $F$). See also \cite{Jez, Diekert_Jez_Pland, EDT0L} for an entirely different approach. Analogous work has been done for other non-solvable groups, such as hyperbolic groups \cite{Rips1995, Dahmani}, partially commutative groups \cite{Casals, Diekert2}, and some free and graph products \cite{Casals2, Diekert}. We refer to \cite{KhMi_icm} for further  results  in this area.  

Note, that there are finitely generated solvable non-virtually abelian groups with decidable Diophantine problem. The first such examples are due to  Kharlampovich, López, and the second author, who  proved in \cite{KLM}  that the Diophantine problem is decidable in the following metabelian groups: $BS(1,n), n \geq 1$ and $A \wr \mbb Z$, where $A$ is a finitely generated abelian group. 

We would like to emphasize that in the case when the Diophantine problem in a group $G$ is undecidable or open it is very interesting to consider decidability of equations or systems of equations of a particular type. In fact, it might be advantageous even in the case when the Diophantine problem in $G$ is decidable, since decision algorithms for particular equations could be much more efficient than the general ones. To this end we would like to mention two results:    it is shown in \cite{quadratic_grigorchuk} that systems of quadratic equations are decidable in the first Grigorchuk group (the Diophantine problem in this group is wide open), and also that orientable quadratic equations are decidable in free metabelian groups \cite{LU} (though the Diophantine problem here is undecidable, see below). \\

\noindent We proceed to state the main results of the paper. In all of them, we consider certain types of groups, and we prove that for any such group $G$ there exists a ring of algebraic integers $O$ that is e-interpretable in $G$, which implies that $\mc{D}(O)$ is polynomial time (many-one) reducible to $\mc{D}(G)$. We further  conjecture that in this case, the ring $\mbb Z$ is e-interpretable in $G$  and  the Diophantine problem in $G$ is undecidable.  In fact, we confirm this conjecture for many  groups $G$ of a particular type.

In Section \ref{subsec:Interpretations} we introduce the notion of interpretation by equations (e-interpretation) of one algebraic structure in another, which is the main technical tool of our method. We show that if a structure $\mathcal{A} $ is e-interpretable in a structure $\mathcal  B$ then $\mc{D}(\mathcal{A})$ is polynomial time many-one reducible to $\mc{D}(\mathcal{B})$, symbolically $\mc{D}(\mathcal{A}) \leq_P \mc{D}(\mathcal{B})$. Reductions of this type are also called \emph{Karp reductions}. All  reductions in this paper are Karp reductions,  so sometimes for brevity we  refer to them simply as  \emph{reductions}.

In Section \ref{subsec:arbitrary nilpotency class} we prove the following principal result,  on which most of the other results of this paper are based on. 
{
\renewcommand{\thetheorem}{\ref{t: main_thm_nilpotent}}
\begin{theorem}\label{t: 1_intro}
Let $G$ be a finitely generated non-virtually abelian nilpotent group. Then there exists a ring of algebraic integers $O$  e-interpretable in $G$, hence $\mc{D}(O)  \leq_P \mc{D}(G)$ (i.e., $\mc{D}(O)$ is Karp reducible to $\mc{D}(G)$).  
\end{theorem}
%\addtocounter{theorem}{-1}
}
As we mention above, if $G$ is virtually abelian then $\mc{D}(G)$ is decidable. 

For nilpotency class $2$, we prove this result by considering the largest ring of scalars $R$ of the bilinear map $G/Z(G)\times G/Z(G) \to G'$ induced by the commutator operation $[\cdot, \cdot]$. It then follows from \cite{GMO_rings} that $R$ is e-interpretable in $G$. By this same reference, there exists a ring of algebraic integers e-interpretable in $R$, and hence in $G$ by transitivity. We refer to Section \ref{s: largest_ring} for further details regarding these results and the notion of largest ring of scalars.  Higher nilpotency class reduces to class $2$ by the following \emph{nilpotent quotient argument}: the third term $\gamma_3(G)$ of the lower central series of  $G$ has finite verbal width, hence it is e-definable in $G$. Therefore the class 2 nilpotent quotient $G/\gamma_3(G)$ is e-interpretable in $G$, which provides the reduction. 

Note that the nilpotent quotient  argument is quite general: it works for any group $G$ with e-definable subgroup $\gamma_3(G)$ or, more generally, $\gamma_i(G), i \geq 3$. 

Theorem \ref{t: 1_intro} together with the nilpotent quotient argument yield the following

{
\renewcommand{\thetheorem}{\ref{t: elliptic_groups}}
\begin{theorem}
Let $G$ be a finitely generated group such that for some $i \in \mathbb{N}$   $\gamma_i(G)$  is e-definable in $G$ (in particular, if $\gamma_i(G)$  has finite verbal width) and  $G/\gamma_i(G)$ is not virtually abelian.  Then there exists a ring of algebraic integers $O$ that is e-interpretable in $G$, hence $\mc{D}(O) \leq_P \mc{D}(G)$.
\end{theorem}
%\addtocounter{\thetheorem}{-1}
}
There are many groups that satisfy the premises of Theorem \ref{t: elliptic_groups}. We mentioned some of them in the corollary below.

{
\renewcommand{\thetheorem}{\ref{c: corollary_quotients}}
\begin{corollary} 
Let $G$ be a finitely generated  group which is either metabelian, or solvable minimax, or polycyclic, or virtually abelian-by-nilpotent, or (nilpotent minimax)-by-(abelian-by-finite). If  $G/\gamma_i(G)$ is not virtually abelian for some $i \in \mathbb{N}$ then there exists a ring of algebraic integers $O$ that is e-interpretable in $G$, hence $\mc{D}(O) \leq_P \mc{D}(G)$.
\end{corollary}
}

The results above are all of the reducibility type, i.e., they show that the Diophantine problem for  various  groups $G$ is at least as hard as the one for  a suitable ring of algebraic integers $O$.  It seems, it does not give us much unless we know the complexity of  the Diophantine problem in $O$, which is conjectured to be undecidable \cite{Phe_Zah, Denef_conjecture}. %Fortunately, 
One can show that for a wide variety of finitely generated nilpotent groups %(and hence the groups $G$ mentioned above) 
the ring $O$ is, in fact, isomorphic to $\mbb Z$, where the Diophantine problem  is known to be undecidable. Furthermore, we proved in \cite{GMO} that the ring $O$ in a random finitely generated nilpotent group of class $c \geq 2$ is, indeed, isomorphic to $\mbb Z$. Our method of proving isomorphism $O \simeq \mbb Z$ is based on maximal rings of scalars of bilinear maps  and  centralizer small   (c-small) elements in groups $G$, i.e., elements of infinite order $g \in G$ such that $C_G(g) = \langle g\rangle \times Z(G)$. 

{
\renewcommand{\thetheorem}{\ref{l: largest_ring_free_2_nilp}}
\begin{proposition}
Let $G$ be a finitely generated nilpotent group of class $2$. If $G$ has a c-small element then the largest ring of scalars of $G$ is $\mbb Z$ and it  is e-interpretable in $G$. 
\end{proposition}
}

We further use the above methods to prove that the ring $\mbb Z$ is e-interpretable in the nilpotent groups below. Let $N_c$ be a  finitely generated non-abelian free nilpotent group of nilpotency class $c \geq 2$.  The next result  is implicit in the aforementioned work of Duchin, Liang and Shapiro \cite{Duchin}. 

{
\renewcommand{\thetheorem}{\ref{t: free_nilpotent}}
\begin{theorem}
The ring $\mbb Z$ is e-interpretable in any finitely generated non-abelian free nilpotent group $N_c$.
\end{theorem}
}

\newcommand{\ut}[1]{UT(#1, \mbb{Z})}
\newcommand{\utt}[2]{UT^{#2}(#1, \mbb{Z})}

{
\renewcommand{\thetheorem}{\ref{t: ut_n_Z}}
\begin{theorem}
Let  $G=\ut{n}$, $n\geq 3$,  be the group of upper uni-triangular $n\times n$ matrices with integer entries. Then the ring $\mbb{Z}$ is e-interpretable in $G$ and $\mc{D}(G)$ is undecidable. 
\end{theorem}
}

We also prove the following
{
\renewcommand{\thetheorem}{\ref{low commutator rank}}
\begin{proposition} Let $G$ be a finitely generated nilpotent group such that $G'/\gamma_3(G)$ has torsion-free rank at most $2$. Then the ring $\mbb Z$ is e-interpretable in $G$, and $\mc{D}(G)$ is undecidable.
\end{proposition} 
%\addtocounter{theorem}{-1}
}
Combining this with a result from \cite{Duchin} we obtain that if $G$ is a finitely generated   non-virtually abelian nilpotent group of class $2$ with infinite cyclic commutator subgroup, then $\mc{D}(G)$ is undecidable, while single equations are decidable in $G$. This applies in particular to any nonabelian generalized Heisenberg group, i.e.\ any group of the form $H_n=\langle a_1,\ldots,a_n,b_1,\ldots,b_n\mid [a_i,b_j]=[a_1,b_1],\ [a_i,a_j]=[b_i,b_j]=1, \ 1\leq i\leq j\leq n\rangle_{\mathcal{N}_2}$ for some $n\geq 1$ (here $\langle \rangle_{\mathcal{N}_2}$ denotes presentation in the variety of nilpotent groups of class $2$). This result was already obtained in \cite{Duchin} for the classical Heisenberg group ($n=1$).

We remark that, to some extent, the converse of Theorem \ref{t: 1_intro} is also true: given a ring of algebraic integers $O$, or more generally, any associative commutative unitary ring $R$, the matrix groups $SL(n,R)$, $T(n,R)$ and $UT(n,R)$ are e-interpretable in the ring $R$ (matrix multiplication can be described by polynomials over $R$). Hence the Diophatine problems in these groups is Karp reducible to the Diophantine problem in $R$.

Finally, we briefly describe the \emph{maximal nilpotent subgroup argument}. In Lemma \ref{max_nilpotent_subgroup} we prove that a maximal finitely generated  nilpotent subgroup $H$ of a fixed class $c$ of a group $G$ is e-definable in $G$. Thus Theorem \ref{t: 1_intro} can also be carried over to groups that have such a (non-virtually abelian) subgroup $H$. With this in mind we denote by ${\mc N}_{max}$ the class of groups where for all $c\geq 1$, every set of  $c$-nilpotent subgroups of $G$ has a maximal element (with respect to inclusion). It turns out that ${\mc N}_{max}$ consists precisely of groups where all abelian subgroups are finitely generated. In Section  \ref{subsec:maximal nilpotent} we describe various types of groups that belong to  ${\mc N}_{max}$, in particular, we note that  groups $GL(n,O)$, where $O$ is a ring of algebraic integers, are there. 

{
\renewcommand{\thetheorem}{\ref{t: maximal_nilp_subgroups}}
\begin{theorem}
Let $ G \in {\mc N}_{max}$.  If  $G$   contains a non-virtually abelian nilpotent subgroup then there exists a  ring of algebraic integers $O$  e-interpretable in $G$, hence $\mc{D}(O) \leq_P \mc{D}(G)$.
\end{theorem} 
}

This result  allows, in particular, to extend Theorem \ref{t: 1_intro} to the class of finitely generated virtually nilpotent groups (since the class ${\mc N}_{max}$ is closed under finite extensions). Another application is the following

{
\renewcommand{\thetheorem}{\ref{c: GL(n,Z)}}
\begin{example}
The Diophantine problem in the groups $GL(3,\mathbb{Z})$, $SL(3,\mathbb{Z})$, and $T(3,\mathbb{Z})$ is undecidable.
\end{example}
}

From Theorem \ref{t: maximal_nilp_subgroups}  and from the fact that any polycyclic group is (nilpotent-by-abelian)-by-finite we obtain the following

{
\renewcommand{\thetheorem}{\ref{t: polycyclic}}
\begin{theorem}
For any virtually polycyclic group  $G$ that is not virtually metabelian there exists a ring of algebraic integers $O$ that is e-interpretable in $G$, and  $\mc{D}(O) \leq_P \mc{D}(G)$.
\end{theorem}
%\addtocounter{\thetheorem}{-1} 
 }
 
Polycyclic metabelian groups constitute an interesting class of groups where to study systems of equations. For instance, if $O$ is a ring of algebraic integers, then $G=O^+\rtimes O^*$ is such a group (by Dirichlet's unit theorem). Here $O^+$ denotes the additive group of $O$ and $O^*$ its group of units, and the action is by ring multiplication. Observe that $G$ is e-interpretable in the ring $O$ (because $O^*$ is e-defined in $O$ by the equation $xy=1$). Hence proving that $\mc{D}(G)$ is undecidable (or that the ring $\mbb Z$ is e-interpretable in $G$) implies that the same is true for the whole ring $O$. An easier task may be to see if  $O$ is e-interpretable in $G$.  Alternatively, it is possible that  $\mc{D}(G)$ is decidable while $\mc{D}(O)$ is  not, e.g.\ if $O=\mbb Z$. Decidability occurs for example in the `trivial' cases when $G$ is virtually abelian (by Lemma \ref{l: Ershov}), which happens if and only if  $O^*$ is finite,  equivalently if $O=\mbb Z$ or $O$ is the ring of integers of an imaginary quadratic field. Recall that a number field is called \emph{quadratic} if it has the form $\mbb Q(\sqrt d)$ for some square-free integer $d$. Such field is said to be \emph{imaginary} if $d<0$, and \emph{real} if $d>0$. Possibly the simplest non-virtually abelian groups of the form $O^+\rtimes O^*$ correspond to real quadratic fields. Hence the next
\begin{problem}\label{problem_intro}
Let $O$ be the ring of integers of a real quadratic number field. Is the Diophantine problem of $O^+\rtimes O^*$ decidable?
\end{problem}

Note that, on the other hand, it is known that $\mc{D}(O)$ is undecidable for $O$ the ring of algebraic integers of any quadratic field \cite{Denef75}.

We finish the paper by studying relatively free groups in the product $\mathcal{A}^d\mathcal{N}_c$ of the varieties $\mathcal{A}^d$ and $\mathcal{N}_c$, where $\mathcal{A}^d$ is the variety of all solvable groups of class $d$ and $\mathcal{N}_c$ is the variety of all nilpotent groups of class $c$, for $c,d\geq 1$. We will refer to these as \emph{free solvable-by-nilpotent groups}. Any such groups is isomorphic to  $F/(\gamma_c(F)^{(d)})$ where $F$ is a free group.

{
\renewcommand{\thetheorem}{\ref{t: solvable_by_nilpotent}}
\begin{theorem}
Let $G$ be a finitely generated nonabelian free (solvable-by-nilpotent) group. Then the ring $\mbb Z$  is e-interpretable in $G$, and $\mc{D}(G)$ is undecidable. 
\end{theorem}
}

The above includes all f.g.\ non-abelian free solvable groups, and it extends Romankov's result that f.g.\ free metabelian groups have undecidable Diophantine problem \cite{Romankov_free_metabelian} (in fact Romankov proves the stronger result that single equations are undecidable in free metabelian groups of countable rank).

{
\renewcommand{\thetheorem}{\ref{t: free_solvable}}
\begin{theorem}
Let $G$ be a finitely generated nonabelian free solvable group. Then the ring $\mbb Z$  is e-interpretable in $G$, and $\mc{D}(G)$ is undecidable. 
\end{theorem}
}

\section{Preliminaries}

\subsection{Interpretations by systems of equations}
\label{subsec:Interpretations}

\paragraph{Multi-sorted structures.} 

A \emph{multi-sorted structure} $\mathcal{A}$ is a tuple $\mathcal{A}=(A_i; f_j, r_k, c_{\ell} \mid i \in I,j \in J,k \in K,\ell\in L),$ where $I,J,K,L$ are sets of natural numbers; the $A_i$ are pairwise  disjoint sets called \emph{sorts}; the $f_j$ are functions of the form
$
f_j: A_{i_{1}} \times \dots \times A_{i_{p}} \to A_{i_p + 1} 
$ for some indices $i_{t}\in I$;
the $r_k$ are relations of the form
$
r_k: A_{s_1} \times \dots \times A_{s_{q}} \to \{0,1\} 
$
for some indices $s_t \in I$; and the $c_{\ell}$ are constants, each one belonging to some sort.  The tuple $(f_j, r_k, c_\ell \mid j,k,\ell)$ is called the \emph{signature} or the \emph{language} of $\mc{A}$. We always assume   that $\mc{A}$ contains the relations "equality in $A_i$" for all sorts $A_i$. 
If $\mc{A}$ has only one sort then $\mc{A}$ is a structure in the usual sense. 
One can construct terms in a multi-sorted structure in an analogous way as in uniquely-sorted structures.  In this case, when introducing a variable $x$, one must specify a sort  where it takes values, which we denote $A_x$.

Let $\mathcal{A}_1, \dots, \mathcal{A}_n$ be a collection of multi-sorted structures.
We let  $(\mathcal{A}_1, \dots,  \mathcal{A}_n)$ be  the multi-sorted structure that is formed by all the sorts, functions, relations, and constants of each $\mathcal{A}_i$. Similarly, given a function $f$ or a relation $r$ we  use the notation $(A,f)$ or $(A,r)$ to refer to the structure $A$ together with the extra function $f$ or relation $r$ in the language, respectively. If two different $\mc{A}_i$'s have the same sort, then we view one of them as a formal disjoint copy of the other.

\paragraph{Diophantine problems and reductions.}\label{Dioph_pblms_intro} 
Let $\mathcal{A}$ be  a multi-sorted structure. An \emph{equation in $\mc{A}$} is an expression  of the form $r(\tau_1, \dots, \tau_k)$, where $r$ is a signature relation of $\mc{A}$ (typically the equality relation), and each $\tau_i$ is a term in $\mc{A}$  where some of its variables may have been substituted by elements of $\mc{A}$. 
Such elements are called the \emph{coefficients} (or the \emph{constants}) of the equation. These may not be signature constants.
A system of equations is a finite conjunction of equations. A \emph{solution} to a system of equations $\wedge_i \Sigma_i(x_1, \dots, x_n)$ on variables $x_1, \dots, x_n$ is a tuple $(a_1, \dots, a_n)\in A_{x_1}\times \dots \times A_{x_{n}}$ such that each $\Sigma_i(a_1, \dots, a_n)$ is true in $\mc{A}$.

The \emph{Diophantine problem} in $\mc{A}$, 
denoted  $\mc{D}(\mc{A})$, refers to the algorithmic problem of  determining if a given system of equations in $\mc{A}$ (with coefficients in $\mc{A}$) has a  solution in  $\mc{A}$.  Sometimes this is  also called \emph{Hilbert's tenth problem} or a \emph{generalized Hilbert's tenth problem} in $\mc{A}$. An algorithm $L$ is  a \emph{decision algorithm} for $\mc{D}(\mc{A})$ if, given a system of equations $S$ in $\mc{A}$, $L$ determines whether or not $S$ has a solution in $\mc{A}$. If such an algorithm exists, then $\mc{D}(\mc{A})$ is called \emph{decidable}, otherwise,  it is \emph{undecidable}. In this paper all structures are finitely generated, the coefficients of equations are given as terms (in the language of $\mc{A}$) in the fixed set of generators.

Let $\mc{A}$ and $\mc{B}$ two structures. We say that Diophantine problem in $\mc{A}$  is \emph{reducible} (or more precisely, \emph{many-one reducible}) to the Diophantine problem in $\mc{B}$ if there is an algorithm that for every finite system of equations $\Sigma$ in $\mc{A}$ constructs a finite system of equations $\Sigma^*$ in $\mc{B}$ such that $\Sigma$ has a solution in $\mc{A}$ if and only if $\Sigma^*$ has a solution in $\mc{B}$.   All our reductions will be polynomial-time computable, i.e., the algorithm that transforms $\Sigma$ to $\Sigma^*$  works in polynomial time. These reductions are called \emph{polynomial-time reductions} or  \emph{Karp reductions}. 
In this case we sometimes write $\mc{D}(\mc{A}) \leq_P \mc{D}(\mc{B})$.

\paragraph{Interpretations by systems of equations.}\label{s: e_interpretations}

Interpretability by systems of equations (e-interpretability) is the analog of the classical notion of interpretability by first-order formulas (see \cite{Hodges, Marker}). In the e-interpretability one requires that only systems of equations are used, instead of arbitrary first-order formulas. As convened above, one is allowed to use any constants (not necessarily from the signature) in such systems of equations.

Let $\mc{A}$ be a structure with sorts $\{A_i \mid i\in I\}$. A \emph{basic set} of $\mc{A}$ is a set of the form
$
A_{i_1} \times \dots \times A_{i_m}
$
for some $m$ and $i_j$'s.

\begin{definition} 
Let $M$ be a basic set of a multi-sorted structure $\mathcal{M}$. A subset $A\subset M$ is called \emph{definable by equations} (or \emph{e-definable}) in $\mathcal{M}$  if there exists a system of equations
$
\Sigma_A(x_1,\ldots,x_m, y_1, \dots, y_k)$ on variables $(x_1, \dots, x_m, y_1, \dots, y_k)=(\mathbf{x}, \mathbf{y})$, such that $\mathbf{x}$ takes values in $M$, and such that 
for any tuple $\mathbf{a}\in M$, one has that  $\mathbf{a} \in A$ if and only if the system $\Sigma_A(\mathbf{a}, \mb{y})$ on variables $\mathbf {y}$ has a solution in $\mathcal{M}$. In this case $\Sigma_A$ is said to \emph{e-define} $A$ in $\mc{M}$.
\end{definition}

From the viewpoint of number theory, an e-definable set is a Diophantine set. From the perspective of algebraic geometry, an e-definable set is a projection  of an affine algebraic set. 
\begin{example}\label{example: center}
Let $G$ be a group generated by $a_1, \dots, a_n$. Then its center $Z(G)$ is e-defined in $G$ by the system of equations $[x, a_i]=1$ ($i=1, \dots, n$) on the variable $x$.
\end{example}

We are ready to introduce the notion of e-interpretability.

\begin{definition}\label{interpDfn}
Let $\mathcal{A}= \left(A_1, \dots; f, \dots, r \dots, c, \dots \right)$  and $\mc{M}$ be two multi-sorted structures. One says that $\mc{A}$ is \emph{interpretable by equations} (or \emph{e-interpretable})  in   $\mathcal{M}$ if for each sort $A_i$ there exists a basic set $M(A_i)$ of $\mathcal{M}$, a subset $X_i\subseteq M(A_i)$,  and an onto map 
$
\phi_i:X_i \to A_i
$ 
such that:
\begin{enumerate}
\item $X_i$ is e-definable in $\mathcal{M}$, for all $i$. 
\item For each function $f$ and each relation $r$ of $\mc{A}$ (including the equality relation of each sort), the preimage by $\boldsymbol{\phi}=(\phi_{1}, \dots)$  of the graph of $f$ (and of $r$) is e-definable in $\mc{M}$, in which case we say that $f$ (or $r$) is e-interpretable in $\mc{M}$. 
\end{enumerate} 
The tuple of maps $\boldsymbol{\phi}=\left(\phi_1, \dots\right)$ is called  an \emph{e-interpretation} of $\mathcal{A}$ in $\mathcal{M}$.% 
\end{definition}

The next two results are fundamental and will often be used without referring to them. They follow from Lemma 2.7 of \cite{GMO_rings}.

\begin{proposition}[\normalfont{E-interpretability is transitive}]\label{interpretation_transitivity}
If $\mc{A}$ is e-interpretable in $\mc{B}$ and $\mc{B}$ is e-interpretable in $\mc{M}$, then $\mc{A}$ is e-interpretable in $\mc{M}$.
\end{proposition}

\begin{proposition}[\normalfont{Reduction of Diophantine problems}]\label{Diophantine_reduction}
Let $\mc{A}$ and $\mc{M}$ be (possibly multi-sorted)  structures such that $\mc{A}$  is e-interpretable in $\mc{M}$.  Then $\mc{D}(O) \leq_P \mc{D}(G)$. As a consequence, if  $\mc{D}(\mc{A})$ is undecidable, then so is $\mc{D}(\mc{M})$.  
\end{proposition}

One of the principal features of e-interpretability is that it is compatible with taking quotients by e-definable congruence relations. Before we see this let us agree on some terminology.

\begin{remark}\label{r: subgroup_edef_eint}
When we say that a subgroup $H$ of a group $G$ is e-definable in $G$, we  mean that $H$ is e-definable as a subset of $G$. Notice that in this case, the identity map $H\to H$ constitutes an e-interpretation of $H$ in $G$. Indeed, the graph of the group operation of $H$ is e-defined in $G$ by the equation $z=xy$, and similarly for the equality relation.
\end{remark}

The following lemma may be read as an illustrative example of the notion of e-interpretability. It can be generalized to any suitable type of structure and the corresponding congrtuence relations.

\begin{lemma}\label{factor_is_interpretable}
Let $N$ be a normal subgroup of a group $G$ such that $N$ is e-definable  in $G$. Then $G/N$ is e-interpretable in $G$. 
\end{lemma}
\begin{proof}
Let $\Sigma(x, \mathbf{y})$ be a system of equations  that e-defines $N$ in $G$, so that $g \in G$ belongs to $N$ if and only if $\Sigma(g, \mathbf{y})$ has a solution $\mathbf{y}$. We check that the natural epimorphism $\pi:G\to G/N$ is an e-interpretation of $G/N$ in $G$. First observe that the preimage of $\pi$ is the whole $G$, which is e-definable in $G$ by an empty system of equations. Regarding equality in $G/N$, the identity  $\pi(g_1) =\pi(g_2)$ holds in $G/N$  if and only if $g_1g_2^{-1}\in N$, i.e.\ if and only if $\Sigma(g_1g_2^{-1}, \mathbf{y})$ has a solution on $\mathbf{y}$. From this it follows that the preimage  of equality in $G/N$,
$
\left\{ g_1, g_2 \in G \mid \pi(g_1) = \pi(g_2) \right\},
$
is e-definable in $G$ by the system $\Sigma(x_1, x_2, \mathbf{y})$ obtained from $\Sigma(x, \mathbf{y})$ after substituting each occurrence of $x$ by $x_1x_2^{-1}$, where $x_1$ and $x_2$ are new variables.  By similar arguments, the preimage of the graph of multiplication in $G/N$ is e-definable in $G$: indeed, $\pi(g_1)\pi(g_2) = \pi(g_3)$ if and only if $g_1g_2g_3^{-1} \in N$.
\end{proof}

\subsection{Varieties of groups}\label{s: varieties}

We will need the following terminology and results. A \emph{class} of groups is a set of isomorphism classes of groups (in this paper any group is identified with its isomorphism class). Given two classes of groups $\mc{S}$ and $\mc{T}$ we let their \emph{product} $\mc{S}\mc{T}$ be the class of all \emph{$\mc{S}$-by-$\mc{T}$} groups, i.e.\ those groups $G$ for which there exists a normal subgroup $N$ such that $N\in \mc{S}$ and $G/N\in \mc{T}$.    A \emph{variety}  of groups $\mc{V}$ is a class of groups for which there exist finitely many  words $w_i(x_1,\dots, x_n)$, $i=1,\dots, m$, on variables $\{x_j\mid j\}$, such that  $K\in \mc{V}$ if and only if $w_i(k_1,\dots, k_n)=1$ for all $k_1, \dots, k_n\in K$ and all $i=1,\dots, m$.  Theorem 21.51 of \cite{Neumann} states that $(\mc{R}\mc{S})\mc{T}=\mc{R}(\mc{S}\mc{T})$ for any three varieties of groups. We let $\mc{A}$, $\mc{E}_n$, $\mc{F}$   denote the classes of all abelian groups, all groups of exponent $n$, and all finite groups,   respectively. The first two are varieties, while the third is not.

\begin{proposition}
Let $H$ be a retract of a group $G$. The the Diophantine problem in $H$ is polynomial time many-one reducible to the Diophantine problem in $G$. In fact the reduction is achieved through the identity map.
\end{proposition}
\begin{proof}
Indeed, let $\pi : G\to H$ be a retract (i.e.\ a homomorphism $G\to H$ which is identical on $H$). Let $S(X)=1$ be an arbitrary finite system of equations with coefficients in $H$, then if $x_i \mapsto a_i\ (x_i\in X)$ is  a solution to this system in $G$, then $x_i\mapsto \pi(a_i)$ is a solution of this system in $H$. Hence $S$ has a solution in $G$ if and only if it has a solution in $H$.
\end{proof}

As an immediate consequence we have the following: 
\begin{corollary}\label{c: countable_uncountable_rank}
Let $\mc{V}$ be a variety of groups and $X$ an infinite set. By $F_{\mc{V}}(X)$ we denote the free group in $\mc{V}$ with basis $X$. Then for any finite subset $X_0\subseteq X$ the subgroup generated by $X_0$ in  $F_{\mc{V}}(X)$ is free in $\mc{V}$ with basis $X_0$, it is a retract of $F_{\mc{V}}(X)$, and its Diophantine problem is polynomial time many-one reducible to the Diophantine problem in $F_{\mc{V}}(X)$.
\end{corollary}

\subsection{Groups and verbal width}\label{s: verbal_width}

As usual we write $[x,y]=x^{-1}y^{-1}xy$ for the \emph{commutator} of two elements $x,y$ of a group $G$, and we let $G^{(0)}=G$, $G^{(1)}=G'=[G,G]=\langle [x,y]\mid x,y\in G \rangle$, and $G^{(i+1)}=(G^{(i)})'$ for all $i\geq 0$. We  further let  $\gamma_1(G)=G$, $\gamma_2(G)=G'$ and $\gamma_{i+1}(G)=[\gamma_i(G), G]$, $i\geq 1$.  The group $G$ is said to be \emph{solvable} of derived length $d\geq 1$ if $G^{(d)}=1$ and $G^{(d-1)}\neq 1$. It is called \emph{nilpotent} of class $c\geq 1$ if $\gamma_{c+1}(G)=1$ and $\gamma_c(G)\neq 1$. An $i$-fold commutator is defined recursively by $[x_1, \dots, x_{i+1}]=[[x_1, \dots, x_i], x_{i+1}]$, for $i\geq 1$. One has $\gamma_i(G)=\langle \{ [g_1, \dots, g_{i}] \mid g_1, \dots, g_{i} \in G \} \rangle$. 

We will need the following auxiliary result.
\begin{lemma}[\hspace{-0.008cm}\normalfont{\cite{Hall}}]\label{l: commutator_finite} 
Any finitely generated finite-by-abelian group  $G$ is  abelian-by-finite (i.e.\ virtually abelian). In fact, in this case the center $Z(G)$ has finite index in $G$.
\end{lemma}

Let $w=w(x_1, \dots, x_m)$ be a word on an alphabet $\{x_1, \dots, x_m\}$. The \emph{$w$-verbal subgroup} of a group $G$ is the subgroup $\langle w(G)\rangle$ generated  by $w(G) = \{ w(g_1,\ldots,g_m)\mid g_1, \dots, g_m\in G \}$. One says that $w$ has \emph{finite width} in  $G$  if there exists an integer $n$ such that every  $g \in \langle w(G)\rangle$ is equal to the product  of at most $n$ elements from $w(G)^{\pm 1}$. In this case  $\langle w(G)\rangle$ is e-defined in $G$ by the equation \begin{equation}\label{e: finite_width}
x=\prod_{i=1}^n \left(w(y_{i1}, \dots, y_{im}) w(z_{i1}, \dots, z_{im})^{-1}\right),
\end{equation}
which has variables $x$ and $\{y_{ij}, z_{ij} \mid 1\leq i\leq n, \ 1\leq j \leq m\}$ (note that some of the factors in \eqref{e: finite_width} can be made trivial by taking $w(1, \dots, 1)^{\pm 1}$). If $w$ has finite width in $G$ for any $w$, then $G$ is said to be \emph{verbally elliptic}. Observe that each term $\gamma_i(G)$ of the lower central series of $G$ is $w_i$-verbal, where $w_i=[x_1,\dots, x_{i}]$.

\begin{remark}\label{r: who_is_verbally_elliptic}
If  $\gamma_i(G)$ has finite width in $G$ then it  is e-definable in $G$  by means of the   equation \eqref{e: finite_width} after taking $[x_1, \dots, x_{i+1}]$ for $w$. Consequently, in this case $G/\gamma_i(G)$ is e-interpretable in $G$, by Lemma \ref{factor_is_interpretable}. 

It is known that any finitely generated nilpotent, metabelian, or polycyclic group is verbally elliptic. More generally, any f.g.\ (abelian-by-nilpotent)-by-finite or (nilpotent minimax)-by-(abelian-by-finite)  group is verbally elliptic. This includes the class of f.g.\ solvable minimax groups. %, and algebraic groups. 
These results are due to George, Romankov, Segal, and Stroud \cite{George, Romankov, Segal2, Stroud}. Proofs can be found in Theorems 2.3.1, 2.6.1, and Corollary 2.6.2 of \cite{Segal2}, respectively. This same reference contains further results of this type for infinitely generated groups. 
\end{remark}

A group $G$ is said to be \emph{minimax} if it admits a composition series all whose factors are finite, infinite cyclic, or quasicyclic (a group is \emph{quasicyclic} if it is isomorphic to $\mbb{Z}[1/p]/\mbb Z$ for some prime $p$). If all the factors are cyclic (finite or infinite), then $G$ is  \emph{polycyclic}.

\section{Largest ring of scalars of  bilinear maps and rings of algebraic integers}\label{s: largest_ring}

Let $A$ and $B$ be abelian groups, and let $f:A\times A\to B$ be a bilinear map between them. We associate with such $f$ a two-sorted structure $(A,B;f)$. The map $f$ is said to be \emph{non-degenerate} if whenever $f(a,x)=0$ for all $x\in A$, one has $a=0$, and similarly for $f(x,a)$. It is  called \emph{full} if the subgroup generated in $B$ by the image of $f$ is  $B$.
An associative commutative unitary ring  $R$  is called a \emph{ring of scalars}  of $f$  if there exist faithful actions of $R$ on $A$ and $B$, which turn $A$ and $B$ into $R$-modules and  such that $f$ is  $R$-bilinear with respect to these actions. More precisely, in this case  $f(\alpha x,y)=f(x,\alpha y)=\alpha f(x,y)$ for all $\alpha\in R $ and all  $x,y\in A$.

Let $R$ be a ring of scalars of $f$. Since $R$ acts faithfully on $A$ and $B$, there exist ring embeddings $R\hookrightarrow \End(A)$ and $R\hookrightarrow \End(B)$. For this reason and for convenience, we  always assume that a ring of scalars of $f$ is a subring of $\End(A)$.  We say that   $R$ is  the \emph{largest} ring of scalars of $f$ if for any other ring of scalars $R'$ of $f$, one has $R'\leq R$ when viewed as subrings of $\End(A)$. If $f$ is full and non-degenerate then such ring exists and is unique \cite{Myasnikov1990}, and  we denote it  $R(f)$. 

The notion of the largest ring of scalars of a bilinear map $f$ was introduced by the second author in \cite{Myasnikov1990}.  This ring constitutes an important feature of $f$, and it has been used successfully to study different first-order theoretic aspects of different types of structures, including rings whose additive group is finitely generated \cite{Mi_So_4}, free algebras \cite{KMLie, KM3, KM2, KM1}, and nilpotent groups \cite{  Mi_so_3, Mi_So_2}.  
For us the most relevant property of $R(f)$ is that it is e-interpretable in $(A, B;f)$:

\begin{theorem}[\normalfont{Theorem 3.5 of \cite{GMO_rings}}]\label{maxringexists}
Let $f:A\times A\to B$ be a full non-degenerate bilinear map between finitely generated abelian groups. Then the largest ring of scalars $R(f)$ of $f$   is finitely generated as an abelian group, and it is e-interpretable in $(A,B;f)$. Moreover $R(f)$ is infinite if and only if $B$ is.
\end{theorem}
\begin{proof}
A more general statement is proved in Theorem 3.5 of \cite{GMO_rings} for $\L$-bilinear maps between $\L$-modules, where $\L$ is an arbitrary  Noetherian   commutative  ring. Our statement corresponds to the  particular case  $\L=\mbb Z$, since the notions of $\mbb Z$-module and of abelian group coincide under the terminology used in \cite{GMO_rings} (see Paragraph 4 of Section 2.3, and  Remark 1.5, both in \cite{GMO_rings}).  
\end{proof}

The previous result constitutes the first step towards e-interpreting rings of algebraic integers in different families of solvable groups (more generally, in structures that have a suitable bilinear map associated to them). The second step is given by the result below. By \emph{rank} of a ring or an abelian group we refer to the maximum number of $\mbb Z$-linearly independent elements in it (considering the group with additive notation). 

\begin{theorem}[\normalfont{Theorem 4.9 and Remark 4.10 of \cite{GMO_rings}}]
Let $R$ be an infinite, finitely generated as an abelian group, associative, commutative unitary ring. Then there exists a ring of algebraic integers $O$ that is e-interpretable in $R$. Moreover, the rank of $O$ does not exceed the rank of $R$.
\end{theorem}

Combining the previous two theorems, we obtain the following fundamental result.

\begin{corollary}\label{c: approach}
Let $f:A\times A \to B$ be a full non-degenerate bilinear map between finitely generated abelian groups, with $B$ infinite. Then there exists a ring of integers $O$ that is e-interpretable in $(A, B; f)$.
\end{corollary}

\section{Diophantine problems in solvable groups}\label{section_applications}

In this section we  present our main results regarding  systems of equations in solvable groups.  
The next lemma deals with the case when the group is virtually abelian.

\begin{lemma}[\normalfont{Proposition 6 of \cite{Ershov}, see also \cite{Noskov}}]\label{l: Ershov}
Any finitely generated virtually abelian group has decidable first-order theory (with constants). In particular, the Diophantine problem in such group is decidable.
\end{lemma}

\subsection{Nilpotent groups}

\subsubsection{Nilpotency class $2$}
In a nilpotent group $G$ of class $2$ the commutator operation $[\cdot, \cdot]$ induces a full non-degenerate bilinear map  between abelian groups:
\begin{equation}\label{e: bilinear_map_of_2_nilp}f: G/Z(G) \times G/Z(G) \to G', \quad \quad  (x Z(G),  yZ(G)) \mapsto [x,y].
\end{equation}
Here $Z(G)$ denotes the center of $G$. By Theorem \ref{maxringexists} the largest ring of scalars $R= R(f)$ of $f$ exists and is e-interpretable in $(G/Z(G), G'; f)$. We denote it by $R(G)$.
\begin{definition}\label{d: largest_ring}
The ring $R = R(G)$ is called the \emph{largest ring of scalars} of $G$.
\end{definition}

Observe that  if $G$ is finitely generated, then both $Z(G)$ and $G'$ are e-definable in $G$ (see Example \ref{example: center} and Remark \ref{r: who_is_verbally_elliptic}, respectively). It follows that the two sorted structure $(G/Z(G), G'; f)$ is e-interpretable in $G$ (indeed the preimage of the graph of $f$ is e-defined in $G$ by the equation $z=[x,y]$). Furthermore if $G$ is not virtually abelian then $G'$ is infinite due to Lemma \ref{l: commutator_finite}. By Theorem \ref{maxringexists}, Corollary \ref{c: approach}, and  transitivity of e-interpretations we obtain the following
\begin{proposition}\label{p: main_prop_2_nilp}
Let $G$ be a finitely generated nilpotent group of nilpotency class $2$. Then the largest ring of scalars $R$ of $G$ is e-interpretable in $G$. If additionally $G$ is not virtually abelian, then there exists a ring of algebraic integers $O$ that is e-interpretable in $R$,  and also in  $G$ by transitivity of e-interpretations. Moreover, the rank of $O$ is at most the rank of $R$.
\end{proposition}

Of particular interest is the case when $G$ is a  finitely generated \emph{free} nilpotent group of nilpotency class $2$.  We shall need the following definition.  

\begin{definition}\label{d: c_small}
An element $g$ in a group $G$ is said to be \emph{c-small} (or centralizer-small)  if $C_G(g)=\{g^t z\mid t\in \mbb Z, \ z\in Z(G)\}$ and $C_G(g)/Z(G)$ is infinite cyclic ($C_G(g)$ denotes the centralizer of $g$ in $G$). 
\end{definition}
We can now prove the following
\begin{proposition}\label{l: largest_ring_free_2_nilp}
Let $G$ be a finitely generated nilpotent group of class $2$. If $G$ has a c-small element then the largest ring of scalars of $G$ is $\mbb Z$ and it is is e-interpretable in $G$. 

\end{proposition}
\begin{proof}
Denote $Z=Z(G)$. Let $a$ be a c-small element of $G$, and let $\psi: G/Z \to G'$ be the group homomorphism given by $xZ \mapsto [a, x]$. Notice that $\psi(xZ)= f(aZ, xZ)$, where $f$ is the map \eqref{e: bilinear_map_of_2_nilp}.
Note also that $\ker(\psi)= C_G(a)/Z = \langle aZ \rangle \cong \mbb Z$. Let $R$ be the largest ring of scalars of $G$. By definition, $R$ acts on $G/Z$ and $G'$, and furthermore $\psi$ is $R$-linear. It follows that $R$ stabilizes $\ker(\psi)\cong \mbb Z$.  Hence for all $\alpha\in R$ there exists an integer $t_\alpha$ such that $\alpha aZ = a^{t_\alpha}Z$. The map $\phi: R \to \mbb Z$  defined by $\alpha \mapsto t_\alpha$ induces a  group embedding between the additive groups of $R$ and $\mbb Z$ (see the proof of Theorem 3.4 in \cite{GMO}). On the other hand, $\mbb Z$ is a ring of scalars of $f$ (note that it acts faithfully on $G/Z$ and $G'$), and so $\mbb Z$ embeds in $R$.  It follows that $R$ as a ring is isomorphic to $\mbb Z$. Finally, Theorem \ref{p: main_prop_2_nilp} implies that $\mbb Z$ is e-interpretable in $G$.
\end{proof}

\subsubsection{Arbitrary nilpotency class}
\label{subsec:arbitrary nilpotency class}

Suppose $G$ is a finitely generated nilpotent group of nilpotency class at least $2$. Then $G/\gamma_3(G)$ is e-interpretable in $G$ and it is nilpotent of class $2$ (see Remark \ref{r: who_is_verbally_elliptic}). Now we can use the methods of the previous section, together with transitivity of e-interpretations, to obtain one of the main results of the paper.

\begin{theorem}\label{t: main_thm_nilpotent}
Let $G$ be a finitely generated non-virtually abelian nilpotent group. Then there exists a ring of algebraic integers $O$  e-interpretable in $G$, and  $\mc{D}(O) \leq_P \mc{D}(G)$. If otherwise $G$ is virtually abelian, then $\mc{D}(G)$ is decidable.
\end{theorem}
\begin{proof}
The last statement of the theorem is  a particular case of Lemma \ref{l: Ershov}.  Hence assume $G$ is not virtually abelian, in which case $G'$ is infinite by Lemma \ref{l: commutator_finite}. This together with Corollary 9 of \cite{Segal} makes $G'/\gamma_3(G)$ infinite as well. It follows that $G/\gamma_3(G)$ is not virtually abelian, since if $G/\gamma_3(G)$ had an abelian subgroup of $G/\gamma_3(G)$ of index say $n$, then $G'/\gamma_3(G)$ would be a finitely generated abelian group of exponent $n^2$, thus finite (to prove this use the identity $[x,y]^k=[x^k,y]=[x,y^k]$, which holds for any two elements $x,y$ in a group of  nilpotency class $2$).   Due to Proposition \ref{p: main_prop_2_nilp}, the largest ring of scalars $R$  of $G/\gamma_3(G)$ is e-interpretable in $G/\gamma_3(G)$, and also in $G$ by transitivity of e-interpretations. Since $G/\gamma_3(G)$ is not virtually abelian, this same proposition implies that  there exists a ring of algebraic integers $O$ that is e-interpretable in $R$, and so also in  $G$.  
\end{proof}

We proceed to study the case of  a finitely generated nonabelian free nilpotent group $N_c$ of arbitrary nilpotency class $c\geq 2$.   The next result  is implicit in the work of Duchin, Liang and Shapiro \cite{Duchin}, and it is made explicit in Corollary 3.3 of \cite{GMO} (which follows the approach of \cite{Duchin}).

\begin{theorem}\label{t: free_nilpotent}
The ring $\mbb Z$ is e-interpretable in any finitely generated non-abelian free nilpotent group $N_c$, hence  $\mc{D}(N_c)$ is undecidable.
\end{theorem}
\begin{proof}
By Proposition \ref{l: largest_ring_free_2_nilp}, the largest ring of scalars of  $N_c/\gamma_3(N_c)$ is $\mbb Z$, and it is e-interpretable in  $N_c/\gamma_3(N_c)$. The latter is e-interpretable in $N$ by Remark \ref{r: who_is_verbally_elliptic}, and hence the result  follows  by transitivity of e-interpretations.  
\end{proof}

The following is an immediate consequence of the above Theorem \ref{t: free_nilpotent} and Corollary \ref{c: countable_uncountable_rank}.
\begin{corollary}
The Diophantine problem is undecidable in any non-abelian free nilpotent group, not necessarily of finite rank.
\end{corollary}

Now we discuss the groups of uni-triangular matrices.
\begin{theorem}\label{t: ut_n_Z}
Let $n\geq 3$ and let $G=\ut{n}$ be the group of upper uni-triangular $n\times n$ matrices with integer entries. Then the ring $\mbb{Z}$ is e-interpretable in $G$ and $\mc{D}(G)$ is undecidable. 
\end{theorem}

\begin{proof}
Given $k>0$, denote by $\utt{n}{m}$ the subgroup of $G$ formed by all matrices of $G$ with $m-1$ zero diagonals above the main one. Then $\utt{n}{1}=\ut{n}$. It is known that for any $r, s > 0$ we have $[\utt{n}{r}, \utt{n}{s}]= \utt{n}{r+s}$ (see Example 3.2.1 in \cite{Merzljakov}). It follows that $\gamma_3(\ut{n})) = \utt{n}{3}$. Let $\overline{G}= G/\gamma_3(\ut{n})$.  

Denote by $t_{ij}$ the transvection matrix with ones on the main diagonal, a one in its $(i,j)$ entry, and  zeros in all other entries. It is straightforward to check that, for all $1\leq i,j,k,\ell \leq n$ with $i\neq j$, $[t_{ik}, t_{k j}]=t_{ij}$ and $[t_{ik}, t_{\ell j}] = 1$  if  $k\neq \ell$. %It follows that $[t_{ik}, t_{\ell i}] = [t_{\ell i}, t_{ik}]^{-1} = t_{\ell, k}^{-1}$ if $k= \ell$ or $[t_{ik}, t_{\ell i}]=1$ if $k\neq \ell$. 
From these identities it follows that  $Z(\overline{G})=\langle \overline{t}_{ij}\mid j-i = 2 \rangle$, $C_{\overline{G}}(\overline{t}_{12}) = \langle \overline{t}_{12}, \overline{t}_{3,4}, \dots, \overline{t}_{n-1,n}, Z(\overline{G})\rangle$, and
$C_{\overline{G}}(\overline{t}_{23}) = \langle \overline{t}_{23}, \overline{t}_{45}, \dots, \overline{t}_{n-1,n}, Z(\overline{G})\rangle$. 

Denote $A=Z(C_{\overline{G}}(\overline{t}_{12}))$ and $B= Z(C_{\overline{G}}(\overline{t}_{23}))$. Again using the above identities  we have that if $n> 5$ then $A=\langle \overline{t}_{12}\rangle Z(\overline{G})$ and  $B=\langle \overline{t}_{23}\rangle Z(\overline{G})$. From $[\overline{t}_{12}, \overline{t}_{23}]=\overline{t}_{13}\in Z(\overline{G})$ we obtain $H=\langle A, B \rangle = A\cdot B$. Hence $H$ is e-definable in $\overline{G}$ and $\overline{t}_{12}$ is a c-small element in $H$. From Proposition  \ref{l: largest_ring_free_2_nilp} and transitivity of e-interpretations we conclude that the ring $\mbb{Z}$ is e-interpretable in $\overline{G}$, and by  Lemma \ref{factor_is_interpretable} and Remark \ref{r: who_is_verbally_elliptic}, we also conclude that $\mbb{Z}$  is  e-interpretable in $G$ for the case $n>5$.

Next we treat the case $3\leq n \leq 5$.  If $n=5$, then using our previous arguments we see that  $A=Z(C_{\overline{G}}(\overline{t}_{12})) = \langle t_{12}\rangle Z(\overline{G})$ and $B=Z(C_{\overline{G}}(\overline{t}_{23})) = \langle t_{23}, t_{45}\rangle Z(\overline{G})$. However, $t_{45}$ belongs to the center of $H=\langle A, B\rangle = A \cdot B$. Therefore $H=\langle t_{12}, t_{23} \rangle Z(\overline{G})$  and $t_{12}$ is a c-small element in $H$. The result then follows similarly as in the case $n> 5$. If $n=4$ then $C_{\overline{G}}(\overline{t}_{23}) = \langle \overline{t}_{23} \rangle Z(\overline{G})$. Hence $\overline{t}_{23}$ is a c-small element of $\overline{G}$, and then the result follows again from Proposition \ref{l: largest_ring_free_2_nilp}.  Finally, the case $n=3$ is already proved in Theorem \ref{t: free_nilpotent}, since then $G$ is free nilpotent.  
\end{proof}

\begin{remark}\label{r: free_nilpotent}

The notion of the largest ring of scalars of a nilpotent group $G$ can be extended to any nilpotency class. This is achieved by considering a bilinear map which resembles the ring multiplication of the  Lie ring of $G$, and which generalizes \eqref{e: bilinear_map_of_2_nilp}. We refer to Subsection 3.3 of \cite{Mi_So_2} for further details, omitting a full explanation here due to its technicality. 

It may be possible to prove that such ring of scalars $R$ is always e-interpretable in $G$.  If so, then the previous results can be approached by considering $R$ directly instead of taking first the quotient $G/\gamma_3(G)$. This approach would yield the same results presented above but with an overall more involved exposition. Nevertheless, it may be more adequate when studying finer aspects of systems of equations in $G$. We shall not pursue this approach in this paper.
 
\end{remark}

We next turn our attention to nilpotent groups $G$ for which $G'/\gamma_3(G)$ has ``small rank''. Recall that by rank of an abelian group or a ring we refer to its maximum number of $\mbb Z$-linearly independent elements (considering the group with additive notation).

\begin{proposition}\label{low commutator rank}
Let $G$ be a finitely generated nilpotent group.  Suppose the rank of $G'/\gamma_3(G)$ is either $1$ or $2$. Then the ring $\mbb Z$ is e-interpretable in $G$,  and $\mathcal{D}(G)$ is undecidable.  
\end{proposition}

\begin{proof}
Let $R$ be the largest ring of scalars of $G/\gamma_3(G)$. By Proposition \ref{p: main_prop_2_nilp}, $R$ is e-interpretable in $G/\gamma_3(G)$, and thus in $G$ as well, by Remark \ref{r: who_is_verbally_elliptic}. By the same arguments as in the proof of Theorem \ref{t: main_thm_nilpotent},  $G/\gamma_3(G)$ is not virtually abelian since otherwise $G'/\gamma_3(G)$ would be finite, contradicting the fact that its rank is $1$ or $2$. Hence  by Proposition \ref{p: main_prop_2_nilp}  there exists a ring of algebraic integers $O$ that is e-interpretable in $R$ with rank at most the rank of $R$. By transitivity, $O$ is e-interpretable in $G$ as well. Let $K$ be the number field of which $O$ is the ring of algebraic integers. It is well known that the rank of $O$ coincides with the degree $|K:\mbb Q|$ of the extension $K/\mbb Q$.   In   \cite{Denef75} Denef  proved that if $K$ is a quadratic field (i.e.\ if $|K:\mbb Q|=2$), then  $\mbb Z$ is e-interpretable in $O$ (see also \cite{Denef_Lipshitz}). Hence if we see that the rank of $R$ is at most $2$, then we will have proved that the ring $\mbb Z$ is e-interpretable in $O$, and also in $G$ by transitivity.

By definition, $R$ acts faithfully by endomorphisms on $G'/\gamma_3(G)$.  
We claim that if the rank of a nontrivial finitely generated abelian group $A$ is at most $2$, then any commutative, associative ring acting faithfully on it also has rank at most $2$. The proof of the proposition will be finished once this claim is proved. 

To prove the claim, first consider the case when $A$ is torsion-free. Then either $A=\mbb Z$ or $A=\mbb Z^2$. If $A=\mbb Z$, then $R\leq \it{End}(A)\cong\mbb Z$ has rank $1$. If $A=\mbb Z^2$, then $R\leq \it{End}(\mbb{Z}^2)$ is a commutative ring whose elements are   $2\times 2$ integer matrices. Let $X$ and $Y$ be two such matrices. Assume that they are not proportional to the identity matrix $I$. Since $X$ and $Y$ commute, elementary calculations show that  $\alpha X+\beta Y + \gamma I = 0$ for some integers $\alpha, \beta, \gamma$ not all of them $0$. This implies that the rank of $R$ is at most $2$ and proves the claim for the case when $A$ is torsion-free.

Now we reduce the general case to the case when $A$ is torsion-free. Let $A$ be a finitely generated abelian group of rank at most $2$, and let $T$  be the torsion subgroup of $A$, i.e.\ the set of all elements of $A$ of finite order. Then $A/T$ is a torsion-free abelian group of rank at most $2$. Notice that for any $a\in T$ and $r\in R$ we have $ra\in T$, hence $R$ acts on $A/T$. Denote $Ann_R(A/T)=\{ r\in R\mid r A \subseteq T\}$. Then $R/Ann_R(A/T)$ acts faithfully on $A/T$, and thus by the paragraph above $R/Ann_R(A/T)$ has rank  at most $2$.  Since $T$ is finite and $A$ is finitely generated, $Hom(A,T)$ is finite (because each homomorphism from $Hom(A,T)$ is uniquely determined by its action on a set of generators of $A$). Hence  $Ann_R(A/T)\leq Hom(A,T)$ is finite. This implies that the rank of $R$ is the same as the rank of $R/Ann_R(A/T)$, which  is at most $2$, and finishes the proof of the claim. 
\end{proof}
%CCCCcheck writing

Recall that the \emph{generalized Heisenberg group} of rank $n$ is the group admitting the following presentation $H_n=\langle a_1,\ldots,a_n,b_1,\ldots,b_n\mid [a_i,b_j]=[a_1,b_1],\ [a_i,a_j]=[b_i,b_j]=1, \ 1\leq i\leq j\leq n\rangle_{\mathcal{N}_2}$, where $\langle \rangle_{\mathcal{N}_2}$ denotes a presentation in the variety of nilpotent groups of class $2$.

\begin{corollary}
Let $G$ be a finitely generated non-virtually abelian nilpotent group of nilpotency class $2$. Assume that $G'$ has rank one. Then the ring $\mbb Z$ is e-interpretable in $G$, and $\mc{D}(G)$ is undecidable. On the other hand,  single equations in $G$ are decidable. This result applies in particular to any generalized Heisenberg group $H_n$ with $n\geq 1$.
\end{corollary}
\begin{proof}
The undecidability part is a particular case of the previous Proposition \ref{low commutator rank}. The decidability part is a consequence of Theorem 3 of \cite{Duchin}. 
\end{proof}

\subsection{Nilpotent quotient argument }

In this section we develop the nilpotent quotient argument mentioned in the introduction and then apply it to Diophantine problems of finitely generated verbally elliptic groups.

\begin{theorem}\label{t: elliptic_groups}
Let $G$ be a finitely generated group such that for some $i \in \mathbb{N}$   $\gamma_i(G)$ is e-definable in $G$ (in particular, if $\gamma_i(G)$  has finite verbal width) and  $G/\gamma_i(G)$ is not virtually abelian.  Then there exists a ring of algebraic integers $O$ that is e-interpretable in $G$, and hence $\mc{D}(O) \leq_P \mc{D}(G)$.
\end{theorem}

\begin{proof}
If for some $i \in \mathbb{N}$   $\gamma_i(G)$ is e-definable in $G$ (for example, if $\gamma_i(G)$ has finite verbal width) and  $G/\gamma_i(G)$ is not virtually abelian,  then    $G/\gamma_i(G)$ is a finitely generated nilpotent not virtually abelian group e-interpretable in $G$.   Hence the conclusions of the theorem hold after applying   Theorem \ref{t: elliptic_groups} and transitivity of e-interpretations.  
\end{proof}

If a group $G$ is verbally elliptic, then $\gamma_i(G)$ has finite verbal width for any natural $i$, and the nilpotent quotient argument applies. All the groups mentioned in the result below are verbally elliptic.

\begin{corollary}\label{c: corollary_quotients}
Let $G$ be a finitely generated  group which is either metabelian, or solvable minimax, or polycyclic, or virtually abelian-by-nilpotent, or (nilpotent minimax)-by-(abelian-by-finite). If  $G/\gamma_i(G)$ is not virtually abelian for some $i \in \mathbb{N}$ then there exists a ring of algebraic integers $O$ that is e-interpretable in $G$, hence $\mc{D}(O) \leq_P \mc{D}(G)$.
\end{corollary}

A definition of a minimax group can be found at the end of Subsection \ref{s: verbal_width}.

\subsection{Groups with maximal nilpotent subgroups}
\label{subsec:maximal nilpotent}

In this section we outline an approach to the Diophantine problem in a group $G$ via the maximal finitely generated nilpotent subgroups of $G$ (if such exist). 
With this in mind we denote by ${\mc N}_{max}$ the class of groups where for all $c\geq 1$, every set of  $c$-nilpotent subgroups of $G$ has a maximal element (with respect to inclusion). Here, and below, for  brevity we  use the expression \emph{$c$-nilpotent} as a replacement of  ``nilpotent of class at most $c$''. 

%. 

\begin{lemma}\label{max_nilpotent_subgroup}
Let $H$ be a finitely generated $c$-nilpotent subgroup of a group $G$, for some $c\geq 1$. Assume  that $H$ is  maximal among all $c$-nilpotent subgroups of $G$. Then $H$ is e-definable  in $G$.
\end{lemma}
\begin{proof}
A group $K$ is $c$-nilpotent if and only if $[k_1,\ldots,k_{c+1}]=1$  for all $k_1, \dots, k_{c+1} \in K$. If $K$ is generated by a finite set, say $e_1, \dots, e_n$, then  this condition holds if and only if $[e_{i_1}, \ldots, e_{i_{c+1}}]=1$ for all $1\leq i_1, \dots, i_{c+1} \leq n$.
We claim that $H$ is e-defined  in $G$ by the following system  of equations on the single variable $x$:
\begin{equation}\label{e: max_nilpotent_def}
\bigwedge_{\substack{z_1, \dots, z_{c+1} \in \\  \left\{x, e_1, \dots, e_m\right\}}} \big( \left[z_1,\ldots,z_{c+1}\right]=1 \big),
\end{equation}
where $\{e_1,\ldots,e_m\}$ is a generating set of $H$. Indeed, since $H$ is $c$-nilpotent, every element  $x \in H$ satisfies \eqref{e: max_nilpotent_def}. Conversely, if $x\in G$ satisfies \eqref{e: max_nilpotent_def} then $\langle H,x\rangle $ is $c$-nilpotent by the observation above. Then by maximality we have $x\in H$. Hence $H$ is e-definable in $G$. 
\end{proof}
The following is the main technical result of our approach to groups from  ${\mc N}_{max}$.

\begin{theorem}\label{t: maximal_nilp_subgroups}
Let $ G \in {\mc N}_{max}$.  If  $G$   contains a non-virtually abelian nilpotent subgroup then there exists a  ring of algebraic integers $O$  e-interpretable in $G$, hence $\mc{D}(O) \leq_P \mc{D}(G)$.
\end{theorem} 

\begin{proof}
Let $H$ be a $c$-nilpotent non-virtually abelian subgroup of $G$.   By Lemma \ref{good_conditions}, $H$ is contained in a subgroup $K\leq G$ that is maximal among all $c$-nilpotent subgroups of $G$, containing $H$. Note that $K$  is also a maximal  $c$-nilpotent subgroup in $G$. This $K$ is  finitely generated due to Condition 3, so it is  e-interpretable in $G$ by Lemma \ref{max_nilpotent_subgroup}. Moreover $K$ is not virtually abelian, since if  $A$ is an abelian  finite index normal subgroup of $K$, then $A\cap H$ is normal in $H$ and abelian, and it has finite index in $H$ (because   $H/ A\cap H$ embeds in $K/A$), a contradiction. The result then follows by Theorem \ref{t: main_thm_nilpotent} and by transitivity of e-interpretations. 
\end{proof}
 
We now provide some properties that guarantee the existence of maximal nilpotent subgroups. 

\begin{lemma}\label{good_conditions}
The following statements are equivalent for any  group $G$:
\begin{enumerate}
\item For all $c\geq 1$, every set of  $c$-nilpotent subgroups of $G$ has a maximal element. 
\item All abelian subgroups of $G$ are finitely generated.
\item All solvable subgroups of $G$ are polycyclic.
\end{enumerate}
\end{lemma}
\begin{proof}
Clearly 1 implies 2, since the set of all finitely generated subgroups of a not finitely generated abelian group has no maximal element.

 Note that if all abelian subgroups of a solvable group are finitely generated, then the group is polycyclic (see, for example, Theorem 21.2.3 from \cite{Kargapolov}). Therefore if Condition 2 holds then every solvable subgroup of $G$ is polycyclic, and so 2 implies 3.

Suppose now Condition 3 holds. Assume that there exists  a set $S$ of $c$-nilpotent subgroups of $G$ with no maximal element, so that $S$ contains an infinite strictly ascending chain of $c$-nilpotent subgroups $N_1<N_2<\dots$ 
The group $N=\bigcup_i N_i$  is $c$-nilpotent, thus it is solvable,  and hence  it is polycyclic, and therefore finitely generated.  This contradicts to the fact that the chain is infinite and strictly ascending.  Hence  3 implies  1.

\end{proof}

\begin{remark}
The result above shows that the class ${\mc N}_{max}$ coincides with the class $A_{max}$ of groups, where each abelian subgroup satisfies $\emph{max}$. The class $A_{max}$ was extensively studied, especially in the case of solvable groups (see \cite{LR}).  
\end{remark}

We next provide some examples of groups from  ${\mc N}_{max}$.  

\begin{example}
Free and torsion free hyperbolic groups are in ${\mc N}_{max}$.
\end{example}

\begin{example}
Limit groups are in ${\mc N}_{max}$ (since centralizers are finitely generated abelian, see \cite{KHARLAMPOVICH1998517}). 
\end{example}

\begin{example}
Any subgroup of $GL(n,O)$ is in ${\mc N}_{max}$  for any $n \in \mathbb{N}$ and any  ring of algebraic integers $O$. 

This is true because any solvable subgroup of $GL(n,O)$ acts faithfully on the additive group of  $O^n$, which is a finitely generated abelian group, and so it is polycyclic (see, for example, Theorem 21.2.2 in \cite{Kargapolov}). 
 \end{example}
 
\begin{example}
Any discrete subgroup of $GL(n,\mathbb{R})$ (with respect to the  topology induced by $\mbb R$) is in ${\mc N}_{max}$.  

Indeed, all discrete solvable subgroups of $GL(n,\mathbb{R})$ are finitely generated  \cite{Auslander}. This is true, in particular, for  abelian subgroups of  a discrete subgroup $G$ of $GL(n,\mathbb{R})$. 

We remark that not all finitely generated subgroups of $GL(n,\mathbb{R})$ are discrete in the induced topology. For instance,  Baumslag-Solitar groups provide examples of finitely generated solvable linear groups that are not polycyclic. 
%. 
\end{example}
%\vspace{0.2cm}

\begin{example}\label{e: example_polycyclic}
Any polycyclic group is in ${\mc N}_{max}$.
\end{example}

\begin{example}\label{e: free_products}
Let $H_1$ and $H_2$ be groups from ${\mc N}_{max}$. Then  $G=H_1*H_2$ is also in  ${\mc N}_{max}$. 

Indeed, by Kurosh subgroup theorem, any subgroup $A$ of $G$ has the form $F*A_1*\cdots *A_m$ for some $m\geq 1$, some free group $F$, and some subgroups   $A_1, A_2$ each of them conjugates to a subgroup of either  $H_1$ or $H_2$.  In particular, if $A$ is abelian then either $A$ is infinite cyclic or $A$ is a conjugate of a subgroup of either  $H_1$ or $H_2$. In both cases $A$ is finitely generated, given that all abelian subgroups in the groups $H_1$ and $H_2$ are finitely generated.

\end{example}
\begin{example}\label{e: direct_products}

Let $G_1, G_2 \in {\mc N}_{max}$. Then $G = G_1 \times G_2$ is in ${\mc N}_{max}$.

Indeed, let  $A$ be an abelian subgroup of $G$. The canonical projections $\pi_1(A)$ and $\pi_2(A)$ on $G_1$ and $G_2$ are abelian. Since $G_1, G_2 \in {\mc N}_{max}$ the subgroups $\pi_1(A)$ and $\pi_2(A)$ are contained in some maximal abelian subgroups $A_1 \leq G_1$ and $A_2 \leq G_2$, which are  finitely generated. It follows that the abelian subgroup $A_1 \times A_2$ of $G$ is also finitely generated. To finish the proof it suffices to note that $A \leq A_1 \times A_2$. 
\end{example}

\begin{example} \label{e:extensions}
Let $G$ be a group and $N$ a normal subgroup of $G$. If $N$ and $G/N$ are in ${\mc N}_{max}$ then $G$ is in ${\mc N}_{max}$.

Indeed, denote by $g \to \bar g $  the canonical epimorphism  $G \to G/N$. Let $A$ be an abelian subgroup of $G$. Then the image $\bar A$ of $A$ in $G/N$ is abelian, and hence finitely generated, say $\bar A  = \langle \bar a_1, \ldots \bar a_n \rangle$, for some elements $a_1, \ldots, a_n \in A$. The subgroup $A \cap N$ is an abelian subgroup of $N$, so it is finitely generated, say $A \cap N = \langle b_1, \ldots,b_m \rangle$. Then $A = \langle a_1, \ldots,a_n,  b_1, \ldots,b_m \rangle$.
\end{example}

\begin{example} \cite{LR} Let $G \in {\mc N}_{max}$. If $A$ is an abelian normal subgroup of $G$ then $G/A \in {\mc N}_{max}$.
\end{example}

We summarize the discussion above in the following results that shows that the class ${\mc N}_{max}$ is rather wide.
\begin{theorem}\label{c: corollary_examples}
The class ${\mc N}_{max}$ satisfies the following conditions: 

\begin{enumerate}
\item  it contains free and torsion free hyperbolic groups.
\item it contains all limits groups.
\item it contains all finite groups and all polycyclic groups.
\item it contains all subgroups of $GL(n,O)$, for any $n \in \mathbb{N}$ and for  any ring $O$ of algebraic integers.

\item it  contains all discrete subgroups of $GL(n,\mbb{R})$.

\item it is closed under direct  products.
\item it is closed under free products.

\item it is closed under  extensions.  
\item it is closed under quotients by abelian normal subgroups.
\end{enumerate}

\end{theorem}

The following  refinement of Theorem \ref{t: main_thm_nilpotent} is an immediate consequence  of  Theorems  \ref{t: maximal_nilp_subgroups} and \ref{c: corollary_examples}.
\begin{corollary}
Let $G$ be a finitely generated virtually nilpotent group that is not virtually abelian. Then there exists a ring of algebraic integers e-interpretable in $G$, and $\mc{D}(O) \leq_P \mc{D}(G)$. If otherwise, $G$ is virtually abelian then $\mc{D}(G)$ is decidable.
\end{corollary}

\begin{example}\label{c: GL(n,Z)}
The Diophantine problem in the groups $GL(3,\mathbb{Z})$, $SL(3,\mathbb{Z})$, and $T(3,\mathbb{Z})$  is undecidable.
\end{example}
\begin{proof}
Let $G$ be one of the groups $GL(3,\mathbb{Z})$, $SL(3,\mathbb{Z})$, or $T(3,\mathbb{Z})$. Then, as mentioned above, $G \in {\mc N}_{max}$. Observe that $G$  contains a finitely generated non-virtually abelian nilpotent subgroup $UT(3,\mathbb{Z})$ of nilpotency class $2$. Suppose we know that $Z(G)UT(3,\mathbb{Z})$ is a maximal $2$-nilpotent subgroup of $G$.
Then it follows by Lemma \ref{max_nilpotent_subgroup} that $Z(G)UT(3,\mathbb{Z})$ is e-definable in $G$. Notice that $Z(G)$ is also e-definable in $G$. Therefore, the quotient group $N = Z(G)UT(3,\mathbb{Z})/Z(G)$ is e-interpretable in $G$. Since $Z(G) \cap UT(3,\mathbb{Z}) = 1$, the group $N$ is isomorphic to $UT(3,\mathbb{Z})$.  Thus, $UT(3,\mathbb{Z})$ is e-interpretable in $G$.  By Theorem \ref{t: ut_n_Z} the ring $\mbb Z$ is e-interpretable in $UT(n,\mathbb{Z})$, hence in $G$ by transitivity of e-interpretations. Hence the Diophantine problem in $G$ is undecidable.

It remains to prove that $N$ is a maximal $2$-nilpotent subgroup of $G$. Assume towards contradiction that this is not the case, and let $x\in G \smallsetminus N$ be such that $N_x = \langle N, x \rangle$ is $2$-nilpotent. Since $N_x$ is infinite, so is $Z(N_x)$, and thus $Z(N_x)$ contains an element $g$ of infinite order. We can write $g=g_0g_1$ for some $g_0 \in Z(UT(3,\mbb{Z}))$ and some $g_1\in Z(G)$. Since all elements of $Z(G)$ have order at most two, we have $g^2 \in Z(UT(3,\mbb{Z}))=\langle t_{1n}(1)\rangle$, where by $t_{ij}(\alpha)$ we denote the transvection matrix with ones on the diagonal, $\alpha$ in the entry $(i,j)$, and zeros everywhere else ($1\leq i,j\leq n$, $\alpha\in\mbb{Z}$). Hence there exists $\alpha\in \mbb{Z}\setminus\{0\}$ such that $t_{1n}(\alpha) \in Z(N_x)$. From the identity $[t_{1n}(\alpha), x]= 1$ we obtain that $x\in T(3,\mbb{Z})$ (all entries below the main diagonal are zero), and the diagonal of $x$ is $(\varepsilon, \delta, \varepsilon)$ for some $\varepsilon, \delta \in \{-1, 1\}$.  Observe that if $G=SL(3,\mbb{Z})$ then we are done, since then $\varepsilon= \delta=1$  and $N_x = N$, a contradiction. 

Next, denote by $y$ the matrix with diagonal $(1,-1,1)$ and with zeros everywhere else. Observe that $xN_x = yN_x$ (this can be seen by multiplying $x$ by transvections $\{t_{ij}\}_{1\leq i<j\leq n}$ and by elements from  $Z(G)$), which implies that $y\in N_x$. We now use $y$ in order  to prove that $N_x$ is not nilpotent. Indeed, let $t^{(0)}$ be the element from $UT(3,\mbb{Z})$ with entries $t_{12}^{(0)} = t_{23}^{(0)}=1$ and $t_{13}^{(0)} = 2$. Define $t^{(i+1)}=[t^i, y] = [t, y, \overset{i}{\dots}, y]$, for all $i\geq 0$. It follows by induction on $i$ that $t_{12}^{(i)} = t_{23}^{(i)}=2^{i}$ and $t_{13}^{(i)} = 2^{2i +1}$ for all $i\geq 0$. Therefore, $N_x$ cannot be $2$-nilpotent, a contradiction.
\end{proof}

\subsection{Polycyclic groups}
We next consider the Diophantine problem in polycyclic groups. 

\begin{theorem}\label{t: polycyclic}
Let $G$ be a virtually  polycyclic group that is not virtually metabelian. Then there exists a ring of algebraic integers $O$ that is e-interpretable in $G$, and $\mc{D}(O) \leq_P \mc{D}(G)$.

\end{theorem}
 
\noindent \emph{Proof of  Theorem \ref{t: polycyclic}\ }
Let $G$ be as in the statement of the theorem, and let $G_0$ be a  finite-index polycyclic normal subgroup of $G$. Any polycyclic group is (nilpotent-by-abelian)-by-finite (see \textsection 2 Theorem 4  of \cite{Segal}). Thus there exists a  chain of subgroups $N\unlhd H\unlhd G_0 \unlhd G$ such that  $G/G_0$ and $G_0/H$ are finite, $H/N$ is abelian, and $N$ is nilpotent. If $N$ is not virtually abelian then there exists a ring of algebraic integers $O$ that is e-interpretable in $G$, by Corollary \ref{c: corollary_examples} and Examples \ref{e: example_polycyclic}, \ref{e:extensions}. 

We claim that if  $N$ is virtually abelian, then $G$  is virtually metabelian. The proof of the theorem will be complete once this claim is proved. We shall use two observations:
\begin{enumerate}
\item Virtually polycyclic groups of finite exponent are finite (this is well known for nilpotent groups, hence it holds for nilpotent-by-abelian, and for their finite extensions).
\item Any normal subgroup or a quotient of a virtually polycyclic group is again virtually polycyclic. 
\item Any finitely generated finite-by-abelian group is  abelian-by-finite, see Lemma \ref{l: commutator_finite}. %(to prove this statement use Theorem 4.25 and Corollary 2 to Theorem 4.21 of \cite{Robinson_finiteness} to reduce it to the case when the group is finitely generated $2$-nilpotent, and then apply our Lemma \ref{l: commutator_finite}). %
\end{enumerate}
Now suppose $N$ has an abelian normal subgroup $A$ such that $A/N$ is finite of order say $n$. Then $H\in (\mc{A}\mc{E}_n)\mc{A}=\mc{A}(\mc{E}_n\mc{A})$ (since the product of varieties is an associative operation ---see Subsection \ref{s: varieties}). 
By Items 1 and 2  we obtain that $H\in \mc{A}(\mc{F}\mc{A})$. Hence $H\in \mc{A}(\mc{A}\mc{F})$. In particular $H\in \mc{A}(\mc{A}\mc{E}_m)$ for some $m$, and by the same reasons as before we have  that $H\in (\mc{A}\mc{A})\mc{F}$. Thus $G_0\in ((\mc{A}\mc{A})\mc{E}_k)\mc{E}_t$ for some $k, t$, and  similarly as before we obtain $G_0\in (\mc{A}\mc{A})(\mc{E}_k\mc{E}_t)\subset(\mc{A}\mc{A})\mc{F}$. Now the same argument yields  $G\in (\mc{A}\mc{A})\mc{F}$. \\

We refer to the introduction for comments regarding the Diophantine problem of polycyclic metabelian groups. See in particular Problem \ref{problem_intro}.

\subsection{Free solvable-by-nilpotent  groups}

We finish the paper by studying systems of equations in finitely generated \emph{free solvable-by-nilpotent} groups. Recall that by this we mean relatively free groups in the product $\mathcal{A}^d\mathcal{N}_c$ of the varieties $\mathcal{A}^d$ and $\mathcal{N}_c$, where $\mathcal{A}^d$ is the variety of all solvable groups of class $d$ and $\mathcal{N}_c$ is the variety of all nilpotent groups of class $\leq c$, for $c,d\geq 1$. These are precisely the groups of the form $F/(\gamma_c(F)^{(d)})$, where $F$ is a free group.

 The following auxiliary lemma is  an immediate consequence of a result due to Malcev \cite{Malcev}.

\begin{lemma}\label{l: FmodN_in_FmodNprime}
Let $N$ be a normal subgroup of a free group $F$ such that $F/N$ is torsion-free. Then $F/N$ is e-interpretable in $F/N'$.
\end{lemma}
\begin{proof}
Let $G=F/N'$. By \cite{Malcev} we have that $C_G(gN')=N/N'$ for any $gN' \in N/N'$. Thus $N/N'$ is e-definable in $G$, and consequently the quotient $F/N$ is e-interpretable in $G$.  
\end{proof}

We start by studying free solvable groups.

\begin{theorem}\label{t: free_solvable}
Let $G$ be a finitely generated  nonabelian free solvable group. Then the ring $\mbb Z$  is e-interpretable in $G$, and $\mc{D}(G)$ is undecidable. 
\end{theorem}
\begin{proof}
Proceed by induction on the derived length  $d$  of $G$. If $d=2$ then $G$ is a f.g.\ free metabelian group. In this case $G$ is verbally elliptic and $G/\gamma_3(G)$ is a finitely generated free $2$-nilpotent group e-interpretable in $G$ (see Remark \ref{r: who_is_verbally_elliptic}). By Theorem  \ref{t: free_nilpotent}, the ring $\mbb Z$ is e-interpretable in  $G/\gamma_3(G)$, and thus in $G$ by transitivity of e-interpretations. 

Now assume  $d\geq 2$. Note that $G=F/F^{(d+1)}$ for some free group $F$. Hence  $G/G^{(d)}=F/F^{(d)}$ is e-interpretable in $G$,
by Lemma \ref{l: FmodN_in_FmodNprime} taking $N= F^{(d)}$.  The quotient $G/G^{(d)}$  is a finitely generated  free solvable group of derived length $d-1$, e-interpretable in $G$. Thus by induction the theorem holds for $G/G^{(d)}$, and then it holds for $G$ as well by transitivity of e-interpretations. 
\end{proof}
%
%qwa

The previous results can be combined to prove the following generalization of Theorems \ref{t: free_nilpotent} and \ref{t: free_solvable}.

\begin{theorem}\label{t: solvable_by_nilpotent}
The ring $\mbb Z$ is e-interpretable in any nonabelian free (solvable-by-nilpotent) group $G$, and $\mc{D}(G)$ is undecidable. 
\end{theorem}
\begin{proof}
We have $G=F/(\gamma_{c}(F)^{(d)})$ for some nonabelian free group $F$ and some integers $c\geq 1$ and $d\geq 1$. We may assume that $c\geq 2$, otherwise $G$ is a  free solvable group and the result follows by the previous Theorem \ref{t: free_solvable}. Proceed by induction on $d$. If $d=1$  then $G$ is free nilpotent and the result is precisely Theorem \ref{t: free_nilpotent}.  Hence suppose $d\geq 2$ and let $N=\gamma_c(F)^{(d-1)}$. Then by  Lemma \ref{l: FmodN_in_FmodNprime}, $F/N$ is e-interpretable in $F/N'=G$. By induction the ring $\mbb Z$ is e-interpretable in $F/\gamma_c(F)^{(d-1)}= F/N$, and the result now follows by transitivity of e-interpretations. 
\end{proof}

The following is a consequence of the above Theorems \ref{t: free_solvable}, \ref{t: solvable_by_nilpotent} and Corollary \ref{c: countable_uncountable_rank}.
\begin{corollary}
The Diophantine problem is undecidable in any non-abelian free solvable group and in any free solvable-by-nilpotent group, not necessarily of finite rank.
\end{corollary}

\bibliographystyle{plain}
\bibliography{bib.bib}   % name your BibTeX data base

\end{document}